\documentclass{IEEEtran}

\usepackage{algorithmic}
\usepackage{graphicx}
\usepackage{textcomp}
\def\BibTeX{{\rm B\kern-.05em{\sc i\kern-.025em b}\kern-.08em
    T\kern-.1667em\lower.7ex\hbox{E}\kern-.125emX}}

\usepackage{umlaute,eurosym,comment}
\usepackage{amsmath,amssymb,amsfonts,bm,upgreek}
\usepackage{float,epsfig,pstricks,pst-node,subfig}
\usepackage{wrapfig}
\usepackage{ifpdf,cite}
\usepackage{multirow,dcolumn}
\usepackage{tikz,pgfplots}

\psset{unit=1mm,framesep=10pt,fillstyle=none}  
\degrees[360]                                  
\psset{linewidth=0.9pt}                        

\newtheorem{theorem}{Theorem}
\newtheorem{lemma}[theorem]{Lemma}

\newtheorem{corollary}[theorem]{Corollary}
\newtheorem{definition}[theorem]{Definition}
\newtheorem{remark}[theorem]{Remark}
\newtheorem{assumption}[theorem]{Assumption}

\providecommand{\ef}{\;.}

\providecommand{\fa}{\forall\;}         

\providecommand{\tp}{^{\scriptsize{\mathrm{T}}}} 

\providecommand{\to}{\rightarrow}                        

\newcommand{\conv}[1]{\text{conv}\left(#1\right)}        

\newcommand{\BZ}{\mathbb{Z}}        
\newcommand{\BR}{\mathbb{R}}        

\newcommand{\BX}{\mathbb{X}} 
\newcommand{\BU}{\mathbb{U}}
\newcommand{\BD}{\mathbb{D}} 
\newcommand{\CX}{\mathcal{X}} 
\newcommand{\CP}{\mathcal{P}} 

\definecolor{ULOcean}{RGB}{0,75,90}
\definecolor{ULOrange}{RGB}{236,116,4}
\definecolor{LiteRed}{RGB}{255,206,206}
\definecolor{MedRed}{RGB}{204,24,24}
\definecolor{DarkRed}{RGB}{127,15,15}
\definecolor{LiteYellow}{RGB}{255,255,100}
\definecolor{MedYellow}{RGB}{255,255,25}
\definecolor{DarkYellow}{RGB}{225,225,0}
\definecolor{LiteBlue}{RGB}{170,184,255}
\definecolor{MedBlue}{RGB}{10,20,204}
\definecolor{DarkBlue}{RGB}{27,29,120}
\definecolor{LiteGreen}{RGB}{24,204,75}
\definecolor{MedGreen}{RGB}{0,130,37}
\definecolor{DarkGreen}{RGB}{0,61,17}
\definecolor{LiteViolet}{RGB}{110,0,225}
\definecolor{MedViolet}{RGB}{73,0,148}
\definecolor{DarkViolet}{RGB}{53,0,107}
\definecolor{LiteOrange}{RGB}{255,226,117}
\definecolor{MedOrange}{RGB}{255,206,0}
\definecolor{DarkOrange}{RGB}{204,167,10}
\definecolor{LiteGray}{RGB}{230,230,230}
\definecolor{MedGray}{RGB}{160,160,160}
\definecolor{DarkGray}{RGB}{100,100,100}

\begin{document}
\title{Deadbeat Robust Model Predictive Control:\newline Robustness without Computing Robust Invariant Sets}
\author{Georg Schildbach, \IEEEmembership{Member, IEEE}
\thanks{Submitted to IEEE TAC on May 11, 2023. This work has been supported by the German Research Foundation (DFG) under grant no.\ 460891204.}
\thanks{Georg Schildbach is with the University of L{\"u}beck, Ratzeburger Allee 160, 23568 L{\"u}beck, Germany (e-mail: georg.schildbach@uni-luebeck.de).}}

\maketitle

\begin{abstract}
Deadbeat Robust Model Predictive Control (DRMPC) is introduced as a new approach of Robust Model Predictive Control (RMPC) for linear systems with additive disturbances. Its main idea is to completely extinguish the effect of the disturbances in the predictions within a small number of time steps, called the deadbeat horizon. To this end, explicit deadbeat input sequences are calculated for the vertices of the disturbance set. They generalize to a nonlinear disturbance feedback policy for all disturbances by means of a barycentric function. Similar to previous approaches, this disturbance feedback policy can be either part of the online optimization (Online DRMPC) or pre-calculated during the design phase of the controller (Offline DRMPC). The main advantage over all other RMPC approaches is that no Robust Positive Invariant (RPI) set has to be calculated, which is often intractable for systems with higher dimensions. Nonetheless, for Online DRMPC and Offline DRMPC recursive feasibility and input-to-state stability can be guaranteed. A small numerical example compares the two versions of DRMPC and demonstrates that the performance of DRMPC is competitive with other state-of-the-art RMPC approaches. Its main advantage is its easy extension to linear time-varying (LTV) and linear parameter-varying (LPV) systems.
\end{abstract}



\vspace*{-0.4cm}
\section{Introduction}\label{Sec:Intro}

Over the past decade, Model Predictive Control (MPC) has developed into one of the most important methods of advanced control. Practical applications have revealed imperfections of nominal MPC regarding different types of uncertainties, such as disturbances, unknown system parameters, and general model mismatch \cite{Mayne:2000}. In many cases, these uncertainties are not negligible for control design. Even though nominal MPC is known to possess a natural level of robustness \cite{ChenAllg:1998,LimonEtAl:2009}, an important focus has developed on Robust MPC (RMPC). RMPC is often able to successfully robustify the controller against these types of uncertainty.

Numerous Robust MPC (RMPC) formulations have therefore been proposed. The most common approach is reformulate the actual system as a linear control system with persistent, additive disturbances inside a bounded set \cite{BemMor:1999,CannonEtAl:2011b}. The disturbance set is constructed by worst-case approximations of the actual uncertainty. The main control task is then to keep the states and the inputs of the linear system inside their respective constraint sets, at all times and for all possible realizations of the uncertainty, which also constitutes the basis for this paper.

\vspace*{-0.2cm}
\subsection*{Literature Review}

Several approaches for solving the RMPC problem have been proposed in the literature \cite{RaMaDi:2018}. Parts of them deal, often quite successfully, with uncertainty in an ad-hoc manner. However, the focus here is only on those with system-theoretic guarantees of recursive feasibility and stability.

One of the most popular approaches fall under Tube-based RMPC (Tube-RMPC) \cite{Mayne:2005,RacEtAl:2012,RaMaDi:2018}. The main idea is to find a linear feedback gain that stabilizes the system and a corresponding minimal robust positive invariant (mRPI) set for the state \cite{RacEtAl:2005}. The control input is split into two parts: one that maintains system state within the mRPI set and one that drives the center of these sets to the origin. The main problems of Tube-RMPC lie in the (usually sub-optimal) selection of the linear feedback gain and the computation (or rather approximation) of the mRPI set.

A second category of RMPC approaches involve a maximal robust positive invariant (MRPI) set. These approaches typically show a slightly lower performance than Tube-RMPC, in terms of the size of the region of attraction. However, the MRPI set has practical advantages over the mRPI set, because exact algorithms are available for its computation and it is tractable even for higher state dimensions \cite{Kerrigan:2000}. A reasonable performance requires closed-loop predictions \cite{Bemporad:1998}. For this to be computationally tractable, piecewise-affine feedback policies are generally used \cite{Loefberg:2003a}, and they should be based on disturbances rather than states (Affine Disturbance Feedback, ADF) \cite{Goulart:2006}.

A third category of approaches is Minimax RMPC \cite{ScoMay:1998}, which is concerned with minimizing the maximal cost over all possible disturbances. Its finite scenario tree is based on the extremal realizations of the disturbances over the prediction horizon. However, this approach suffers from exponential complexity with growing prediction horizon \cite{Loefberg:2003a}.

\vspace*{-0.5cm}
\subsection*{Contributions}

Many refinements of the above RMPC approaches have been proposed in the last decade. In this paper, a completely new approach is explored, called Deadbeat Robust Model Predictive Control (DRMPC). It bears some similarities, but also some differences with each of the three categories. The rationale is that DRMPC offers benefits for some particular applications, and may serve as a useful basis for control problems that remain challenging today, like linear time-varying (LTV) and linear parameter-varying (LPV), but also nonlinear and large-scale systems.

Coming from Minimax MPC \cite{ScoMay:1998,Loefberg:2003a}, DRMPC introduces a second horizon, significantly smaller than the prediction horizon, over which the effect of the disturbances is completely extinguished (see Section \ref{Sec:Background}). For a controllable linear system, this so-called \emph{deadbeat horizon} can be in the order of the number of states, so that the combinatorial growth of the scenario tree remains manageable. In contrast to Minimax MPC approaches, an exact solution may hence be possible.

Compared to ADF, DRMPC is not restricted to affine feedback policies in the predictions. Instead, it features a full nonlinear feedback policy based on a barycentric function. In terms of the feasibility of linear MPC with convex constraints, this is the most general class of nonlinear policies that one can wish for. Similar to ADF, DRMPC allows to solve for the optimal feedback policies online (`Online DRMPC'), while retaining convexity of the optimal control problem (see Section \ref{Sec:DRMPC}). Yet ADF requires even a Robust Control Invariant (RCI) set for this purpose \cite{Goulart:2006}. Alternatively, the feedback policies in DRMPC can also be pre-computed offline (`Offline DRMPC'). Then the solution is sub-optimal, but it reduces the computational complexity of DRMPC to approximately the same level as nominal MPC. The Online DRMPC and Offline DRMPC approaches both come with full guarantees of recursive feasibility and stability (see Section \ref{Sec:Guarantees}).

Similar to Tube-RMPC, DRMPC also splits the control input into two parts: one that compensates the disturbances and one that drives the state to the origin on a bigger scale. Unlike Tube-RMPC, however, DRMPC does not require the a priori selection of a stabilizing feedback gain. In fact, a linear feedback gain is not required at all. Instead, the feedback policy and hence the exact split of the control input is computed optimally. Furthermore, DRMPC does not require the a priori computation of an mRPI set for the tube section. For the terminal set, a (non-robust) positive invariant (PI) set suffices, e.g., simply the origin. This is of practical relevance because the origin is a popular choice especially for large scale systems, where (robust) positive invariant sets are too complex or computationally intractable. Even more importantly, this setup makes the extension of DRMPC to the robust control of LTV and LPV systems straightforward. 

\vspace*{-0.2cm}
\section{Background and Terminology}\label{Sec:Background}


\vspace*{-0.1cm}
\subsection{Basic Terminology and Notation}\label{Sec:Notation}

$\mathbb{R}$ and $\mathbb{Z}$ denote the sets of real and integral numbers. $\mathbb{R}_{+}$ ($\mathbb{R}_{0+}$) and $\mathbb{Z}_{+}$ ($\mathbb{Z}_{0+}$) are then the sets of positive (non-negative) real and integral numbers, respectively. 

A function $\alpha:\BR_{0+}\rightarrow\BR_{0+}$ is a \textbf{$\textbf{K}$-function} if it is continuous, strictly monotonically increasing, and $\alpha(0)=0$. It is a \textbf{$\textbf{K}_{\infty}$-function} if, in addition, $\alpha(r)\rightarrow\infty$ as $r\rightarrow\infty$. A function $\beta:\BR_{0+}\times\BZ_{0+}\rightarrow\BR_{0+}$ is a \textbf{$\textbf{KL}$-function} if $\beta(\,\cdot\,,k)$ is a $\text{K}$-function for any fixed $k\in\BZ_{0+}$, and $\beta(r,\,\cdot\,)$ is monotonically decreasing with $\beta(r,k)\rightarrow 0$ as $k\rightarrow\infty$ for any fixed $r\in\BR_{0+}$.

On a vector space $\mathbb{R}^{n}$, for any $n\in\mathbb{Z}_{+}$, $\|\cdot\|$ is used to denote any norm. If $\mathbb{S},\mathbb{T}\subset\mathbb{R}^{n}$ are sets, the \textbf{Minkowski sum} is defined as
\begin{equation*}
	\mathbb{S}\oplus \mathbb{T}\triangleq\bigl\{s+t\in\mathbb{R}^{n}\:\big|\: s\in \mathbb{S},\:t\in \mathbb{T}\}\;,
\end{equation*}
and the \textbf{Pontryagin difference} as
\begin{equation*}
	\mathbb{S}\ominus \mathbb{T}\triangleq\bigl\{\xi\in\mathbb{R}^{n}\:\big|\:\xi+t\in \mathbb{S}\:\:\forall\:t\in \mathbb{T}\}\;.
\end{equation*}

\vspace*{-0.4cm}
\subsection{Control System}\label{Sec:System}

Consider the discrete-time (DT), linear-time invariant (LTI) control system with additive disturbances
\begin{equation}\label{Equ:DTSystem}
  x_{k+1}=Ax_{k}+Bu_{k}+d_{k}\;.
\end{equation}
Here $k\in\mathbb{Z}_{0+}$ denotes the time step, $x_{k}\in\mathbb{R}^{n}$ and $u_{k}\in\mathbb{R}^{m}$ are the state and the input and $d_{k}\in\mathbb{R}^{n}$ is the disturbance at time step $k$, respectively. System \eqref{Equ:DTSystem} comes with a given initial condition $x_{0}\in\mathbb{R}^{n}$.

\begin{assumption}[control system]\label{The:System}
  (a) The state of the system $x_{k}$ can be measured in each step $k$. (b) The system matrix $A\in\BR^{n\times n}$ and input matrix $B\in\BR^{n\times m}$ form a controllable pair.
\end{assumption}

The control objective is to regulate the state of system \eqref{Equ:DTSystem} to the origin, while respecting state and input constraints:
\begin{equation}\label{Equ:Constraints}
  x_{k}\in\BX\subseteq\BR^{n}\;,\enspace u_{k}\in\BU\subseteq\BR^{m}\quad\fa k\in\mathbb{Z}_{0+}\ef
\end{equation}
The additive disturbance is also contained in a given disturbance set:
\begin{equation}\label{Equ:DisSet}
  d_{k}\in\BD\subseteq\BR^{n}\quad\fa k\in\mathbb{Z}_{0+}\;.
\end{equation}

\begin{assumption}[constraint sets]\label{The:Constraints}
  (a) The input and state constraint sets $\BU$ and $\BX$ are convex and contain the origin. (b) The disturbance set $\BD$ is a convex polytope that contains the origin.
\end{assumption}

\vspace*{-0.4cm}
\subsection{Deadbeat Input Sequence}\label{Sec:DeadbeatInputSequence}

This section gradually introduces the concept of an \emph{deadbeat input sequence}, which is the main idea behind DRMPC.

For a start, consider the control problem without constraints. Suppose that the state of the system in some time step $k-1\in\BZ_{0+}$ is at the origin, $x_{k-1}=0$, whence $u_{k-1}=0$. Consider only the effect of a single disturbance in step $k-1$, $d_{k-1}\neq 0$. According to \eqref{Equ:DTSystem}, the state of the system at step $k$ is given by $x_{k}=d_{k-1}$.

By Assumption \ref{The:System}(b), the system state can be steered from $x_{k}$ back to the origin in at most $n$ steps, assuming there no further disturbances during these steps. So pick any $M\geq n$, and denote with $z_{k|k},\hdots,z_{k+M|k}$ the \emph{predicted states} under the assumption of $d_{k}=\hdots=d_{k+M-1}=0$. There hence exists a sequence of inputs $w_{k|k},\hdots,w_{k+M-1|k}$ such that
\begin{equation*}
  z_{k+j+1|k}=Az_{k+j|k}+Bw_{k+j|k}\quad\fa j=0,\hdots,M-1\;,
\end{equation*}
with $z_{k|k}=x_{k}=d_{k-1}$ and $z_{k+M|k}=0$. 

\begin{definition}[deadbeat input sequence]
  (a) The input sequence $w_{k|k},\hdots,w_{k+M-1|k}\in\BR^{m}$ is called an \textbf{deadbeat input sequence} for the disturbance $d_{k-1}$. (b) The sequence $z_{k|k},\hdots,z_{k+M|k}\in\BR^{n}$ is called the corresponding \textbf{state predictions}. (c) $M\in\BZ_{+}$ is called the \textbf{deadbeat horizon} ($M\geq n$).
\end{definition}

Next, suppose there is another disturbance in step $k$, $d_{k-1}\neq 0$ and $d_{k}\neq 0$. Let $w_{k+1|k+1},\hdots,w_{k+M|k+1}$ be a deadbeat input sequence for $d_{k}$. By the superposition principle for linear systems, the following selection of inputs

\begin{align*}
  u_{k}&=w_{k|k}\;,\\
  u_{k+1}&=w_{k+1|k}+w_{k+1|k+1}\;,\\
  &\;\;\vdots\\
  u_{k+M-1}&=w_{k+M-1|k}+w_{k+M-1|k+1}\;,\\
  u_{k+M}&=w_{k+M|k+1}
\end{align*}
annihilates both disturbances $d_{k-1}$ and $d_{k}$, in the sense that the predicted state $z_{k+M+1|k+1}$ is zero.

This procedure can be generalized to the case where all disturbances are non-zero, i.e., $d_{k}\neq 0$ for all $k\in\BZ_{0+}$. By the considerations above, the required input to be applied at step $k$ is given by
\vspace*{-0.2cm}
\begin{equation}\label{Equ:AnnihilatingInput}
  u_{k}=\sum_{j=0}^{\max\{k-1,M-1\}}w_{k|k-j}\;,
\end{equation}
\vspace*{-0.1cm}
where $w_{k|k-j}$ is the element of the deadbeat input sequence for the disturbance $d_{k-j-1}$ to be applied at step $k$.

\vspace*{-0.4cm}
\subsection{Vertex Representations}\label{Sec:VertexRepresentations}

Assumption \ref{The:Constraints}(b) implies that the disturbance set can be represented as the \emph{convex hull} of a finite number $i=1,\hdots,p$ of vertices $d^{(i)}$,
\vspace*{-0.1cm}
\begin{equation}\label{Equ:ConvHull}
  \BD=\conv{d^{(1)},d^{(2)},\hdots,d^{(p)}}\;.
\end{equation}
\vspace*{-0.5cm}

\noindent For any $k\in\BZ_{0+}$, the disturbance vector $d_{k}\in\BD$ can hence be expressed as
\vspace*{-0.1cm}
\begin{equation}\label{Equ:DisConvHull}
  d_{k}=\sum_{i=1}^{p}\lambda_{k}^{(i)}d^{(i)}\;,
\end{equation}
\vspace*{-0.3cm}

\noindent for appropriate scalars $\lambda_{k}^{(i)}\geq 0$ with $\sum_{i=1}^{p}\lambda_{k}^{(i)}=1$.

Denote with $w_{0}^{(i)},\hdots,w_{M-1}^{(i)}$ the deadbeat input sequence for the case where $x_{0}$ coincides with the vertex point $d^{(i)}$. Due to linearity of the system, the deadbeat input sequence for \emph{any} $d_{k-1}\in\BD$ can be expressed as
\vspace*{-0.1cm}
\begin{equation}\label{Equ:InputConvHull}
  w_{k+j|k}=\sum_{i=1}^{p}\lambda_{k-1}^{(i)}w^{(i)}_{j}\quad\fa j=0,\hdots,M-1\;,
\end{equation}
\vspace*{-0.3cm}

\noindent where $\lambda_{k-1}^{(i)}$ are the scalars in \eqref{Equ:DisConvHull} pertaining to $d_{k-1}$.

Moreover, denote with $z_{0}^{(i)},\hdots,z_{M}^{(i)}$ the corresponding state predictions for the vertex points $x_{0}=d^{(i)}$:
\begin{equation*}
  z_{j+1}^{(i)}=Az_{j}^{(i)}+Bw_{j}^{(i)}\;,\quad\fa j=0,\hdots,M-1\;.
\end{equation*}
Note that $z_{0}^{(i)}=d^{(i)}$ and $z_{M}^{(i)}=0$ for all $i=1,\hdots,p$ due to the deadbeat property of the input sequences. By the superposition principle, the state predictions $z_{k|k},\hdots,z_{k+M|k}$ corresponding to \eqref{Equ:InputConvHull} can be expressed as
\begin{equation}\label{Equ:StateConvHull}
  z_{k+j|k}=\sum_{i=1}^{p}\lambda_{k-1}^{(i)}z^{(i)}_{j}\quad\fa j=0,\hdots,M\;,
\end{equation}
\vspace*{-0.3cm}

\noindent where $\lambda_{k-1}^{(i)}$ are the scalars in \eqref{Equ:DisConvHull} pertaining to $d_{k-1}$. Observe that $z_{k|k}=d_{k-1}$ and $z_{k+M|k}=0$ by construction, because $z_{M}^{(i)}=0$ for all $i=1,\hdots,p$.

\begin{remark}\label{Rem:NonlinearDF}
  The inputs $w_{k|k},\dots,w_{k+M-1|k}$ in \eqref{Equ:InputConvHull} are given as nonlinear functions of $d_{k-1}$ (`\textbf{disturbance feedback policy}'), defined through the barycentric coordinates $\lambda_{k-1}^{(1)},\dots,\lambda_{k-1}^{(p)}$ of the vertices.
\end{remark}

\vspace*{-0.1cm}
\section{Deadbeat Robust MPC}\label{Sec:DRMPC}

\subsection{Constraint Tightening}\label{Sec:Constraints}

For simplicity, consider the DRMPC problem setup only in step $k=0$. The generalization to an arbitrary step $k=1,2,\dots$ is straightforward. Let $x_{0}\in\BX$ be the initial condition of the system and $N$ be the selected \emph{prediction horizon}, where $n\leq M\leq N$. By Assumption \ref{The:System}(b), there exist control inputs $u_{0|0},\hdots,u_{N-1|0}$ that drive the system state to the origin; i.e., the predicted states $x_{0|0},\hdots,x_{N|0}$, calculated recursively via\vspace*{-0.1cm}
\begin{equation*}
  x_{k+1|0}=Ax_{k|0}+Bu_{k|0}\;,\quad x_{0|0}=x_{0}\;,
\end{equation*}
\vspace*{-0.5cm}

\noindent satisfy $x_{N|0}=0$. Thus the terminal set for the DRMPC is chosen as the origin.

The notation $u_{k+j|k}$ refers to the input to be applied to the system at step $k+j$, according to the plan made in step $k$. Correspondingly, $x_{k+j|k}$ refers to the state prediction for step $k+j$ made in step $k$.

According to \eqref{Equ:AnnihilatingInput}, for any disturbance $d_{k-1}$ that is observed at any step $k=1,\hdots,N-1$, the plan is to superimpose a deadbeat input sequence $w_{k|k},\hdots,w_{k+M-1|k}$ over the subsequent $M$ inputs, in order to eliminate the effects of the disturbance. For instance,\vspace*{-0.2cm}
\begin{subequations}\label{Equ:DistCompensation}\begin{align}
  u_{0}&=u_{0|0}\;,\\
  u_{1}&=u_{1|0}\,\,+\,\,w_{1|1}\;,\\
  u_{2}&=u_{2|0}\,\,+\,\,w_{2|1}\,\,+\,\,w_{2|2}\;,\\
  &\;\;\vdots\nonumber\\
  u_{M}&=u_{M|0}+w_{M|1}+\hdots+w_{M|M}\;,\hspace*{0.9cm}\phantom{a}
\end{align}\end{subequations}
\vspace*{-0.5cm}

\noindent up to the deadbeat horizon, and\vspace*{-0.1cm}
\addtocounter{equation}{-1}
\begin{subequations}\setcounter{equation}{4}\begin{equation}
  u_{k}=u_{k|0}+\sum_{j=0}^{M-1}w_{k|k-j}\quad\fa M\leq k<N\;.
\end{equation}\end{subequations}
\vspace*{-0.3cm}


Since the disturbances $d_{0},\hdots d_{N-1}$ over the prediction horizon are unknown at step $0$, the disturbance compensation terms $w_{k|k-j}$ in \eqref{Equ:DistCompensation} are not fixed when the inputs $u_{0|0},\hdots,u_{N-1|0}$ have to be computed. It is thus necessary to tighten the input constraints accordingly. With $w^{(1)}_{j},\hdots,w^{(p)}_{j}$ from \eqref{Equ:InputConvHull}, the tightening can be based on the corner cases\vspace*{-0.1cm}
\begin{multline}
  w_{k|k-j}\in\conv{w^{(1)}_{j},\hdots,w^{(p)}_{j}}\\ \fa j=0,\hdots,M-1\;.
\end{multline}
\vspace*{-0.5cm}

\noindent The corner case inputs are included as decision variables into the DRMPC optimization problem. The tightened input constraints can hence be written as\vspace*{-0.1cm}
\begin{subequations}\label{Equ:InputTightening}\begin{align}
  &u_{0|0}\in\BU\;,\\
  &u_{1|0}+w^{(i_0)}_{0}\in\BU\;,\\
  &u_{2|0}+w^{(i_0)}_{0}+w^{(i_1)}_{1}\in\BU\;,\\
  &\;\;\vdots\nonumber\\
  &u_{M|0}+w^{(i_0)}_{0}+\hdots+w^{(i_{M-1})}_{M-1}\in\BU\hspace*{2.05cm}\phantom{a}
\end{align}\end{subequations}
\vspace*{-0.5cm}

\noindent up to the deadbeat horizon, and beyond that\vspace*{-0.1cm}
\addtocounter{equation}{-1}
\begin{subequations}\setcounter{equation}{4}\begin{equation}
  u_{k|0}+w^{(i_0)}_{0}+\hdots+w^{(i_{M-1})}_{M-1}\in\BU\quad\fa M\leq k<N\;.
\end{equation}\end{subequations}
The tightening is such that \eqref{Equ:InputTightening} must hold for all combinations of $i_0,\hdots,i_{M-1}\in\{1,\hdots,p\}$. Thus (\ref{Equ:InputTightening}b) represents $p$ input constraints, (\ref{Equ:InputTightening}c) represents $p^2$ input constraints, and (\ref{Equ:InputTightening}d,e) represent $p^M$ input constraints. This sums up to a total of $1+p+p^2+\hdots+p^{(M-1)}+(N-M)p^M$ input constraints.

The tightening of the corresponding state constraints is analogous. Based on the inputs \eqref{Equ:DistCompensation} and the predicted states $x_{0|0},\hdots,x_{N|0}$, the true states satisfy
\begin{subequations}\label{Equ:DistCompensation}\begin{align}
  x_{1}&=x_{1|0}+z_{1|1}\;,\\
  x_{2}&=x_{2|0}+z_{2|1}+z_{2|2}\;,\\
  &\;\;\vdots\nonumber\\
  x_{M}&=x_{M|0}+z_{M|1}+\hdots+z_{M|M}\;,\hspace*{1cm}\phantom{a}
\end{align}\end{subequations}
up to the deadbeat horizon, and beyond that
\addtocounter{equation}{-1}
\begin{subequations}\setcounter{equation}{4}\begin{equation}
  x_{k}=x_{k|0}+\sum_{j=0}^{M-1}z_{k|k-j}\quad\fa M\leq k<N\;.
\end{equation}\end{subequations}

Again, the tightening of the state constraints can be based on the corner cases of the deadbeat state predictions \eqref{Equ:StateConvHull},\vspace*{-0.1cm}
\begin{multline}
  z_{k|k-j}\in\conv{z^{(1)}_{j},\hdots,z^{(p)}_{j}}\\ \fa j=0,\hdots,M-1\;.
\end{multline}
\vspace*{-0.5cm}

\noindent The corner case states are also decision variables of the DRMPC optimization problem. The tightened state constraints can hence be written as
\begin{subequations}\label{Equ:StateTightening}\begin{align}
  &x_{1|0}+z^{(i_0)}_{0}\in\BX\;,\\
  &x_{2|0}+z^{(i_0)}_{0}+z^{(i_1)}_{1}\in\BX\;,\\
  &\;\;\vdots\nonumber\\
  &x_{M|0}+z^{(i_0)}_{0}+\hdots+z^{(i_{M-1})}_{M-1}\in\BX\hspace*{2.05cm}\phantom{a}
\end{align}\end{subequations}
up to the deadbeat horizon, and beyond that\vspace*{-0.1cm}
\addtocounter{equation}{-1}
\begin{subequations}\setcounter{equation}{4}\begin{equation}
  x_{k|0}+z^{(i_0)}_{0}+\hdots+z^{(i_{M-1})}_{M-1}\in\BX\quad\fa M\leq k\leq N\;.
\end{equation}\end{subequations}
\vspace*{-0.5cm}

\noindent The tightening is such that \eqref{Equ:StateTightening} must hold for all combinations of $i_0,\hdots,i_{M-1}\in\{1,\hdots,p\}$. Thus, in total, there are $p+p^2+\hdots+p^{(M-1)}+(N-M+1)p^M$ state constraints, not counting the additional terminal constraint.

\vspace*{-0.1cm}
\subsection{DRMPC Optimization Problem}\label{Sec:DARMPCProblem}

All required constraints have already been derived. The standard cost function consists of a sum of \emph{stage costs} $\ell:\BX\times\BU\to\BR$ over the prediction horizon. Thus the \emph{DRMPC (optimization) problem} reads as follows:
\begin{subequations}\label{Equ:DRMPCProblem}\begin{align}
  \min\enspace&\sum_{k=0}^{N-1}\ell\left(u_{k|0},x_{k|0}\right)\\
  \text{s.t.}\enspace&\text{for all}\enspace j\in\{0,\hdots,M-1\}\;,\;\;i\in\{1,\hdots,p\}:\nonumber\\
   &\quad z^{(i)}_{j+1}=Az^{(i)}_{j}+Bw^{(i)}_{j}\;,\enspace z^{(i)}_{0}=d^{(i)}\;,\\
   &\quad z^{(i)}_{M}=0\;,\\
   &\text{for all}\enspace k\in\{0,\hdots,N-1\}:\nonumber\\
   &\quad x_{k+1|0}=Ax_{k|0}+Bu_{k|0}\;,\enspace x_{0|0}=\bar{x}_{0}\\
   &\quad x_{N|0}=0\\
   &\quad u_{k|0}\in\BU\\
   &\quad x_{k+1|0}\in\BX
\end{align}\end{subequations}

\addtocounter{equation}{-1}

\vspace*{-0.9cm}
\begin{subequations}\begin{align}\setcounter{equation}{7}
   &\quad\text{for all}\enspace i_{0},\hdots,i_{M-1}\in\{1,\hdots,p\}:\nonumber\\
   &\qquad u_{k|0}+\sum_{l=0}^{\max\{k-1,M-1\}} w^{(i_{l})}_{l}\in\BU\;,\\
   &\qquad x_{k+1|0}+\sum_{l=0}^{\max\{k,M-1\}} z^{(i_{l})}_{l}\in\BX\;.
\end{align}\end{subequations}
Recall that the sums in (\ref{Equ:DRMPCProblem}h,i) refer to all possible combinations of the indices $i_{l}$. 

\begin{remark}[terminal condition]\label{Rem:TerminalCondition}
  For the DRMPC formulation \eqref{Equ:DRMPCProblem}, the \textbf{terminal set} is the origin and there is no \textbf{terminal cost}.
\end{remark}

It is possible, however, to choose a different terminal condition for a DRMPC in the case with offline pre-computations, which is the subject of Section \ref{Sec:OfflineDRMPC}.

\vspace*{-0.4cm}
\subsection{Computational Complexity}\label{Sec:DARMPCComplexity}

If the stage cost is a convex function and $\BU$, $\BX$ are convex sets, then \eqref{Equ:DRMPCProblem} is a convex optimization problem. In this case, it can be solved quickly and efficiently, even for larger problem instances \cite{BoydVan:2004,WangBoyd:2010}. There are many excellent software tools available, some of which also support code generation for embedded systems \cite{MattBoyd:2012,Ferreau:2017,ZanelliEtAl:2020,SchallerEtAl:2022}.

The decision variables of the DRMPC problem comprise 
\begin{itemize}
  \item $Nm$ selected inputs $u_{0|0},\hdots,u_{N-1|0}$\;,
  \item $Nn$ predicted states $x_{1|0},\hdots,x_{N|0}$\;,
  \item $Mmp$ deadbeat inputs $w_{0}^{(i)},\hdots,w_{M-1}^{(i)}$\;,
  \item $Mnp$ corresponding states $z_{1}^{(i)},\hdots,z_{M}^{(i)}$\;.
\end{itemize}
The total of $(Mp+N)(m+n)$ decision variables can be reduced to $(Mp+N)m$, if all state variables are eliminated by substitution of the equality constraints (\ref{Equ:DRMPCProblem}b,c) and (\ref{Equ:DRMPCProblem}d,e). This is sometimes called the \emph{sequential approach}.

The total number of inequality constraints in \eqref{Equ:DRMPCProblem} is\vspace*{-0.1cm}
\begin{multline}\label{Equ:IneqConstraints}
  \hspace*{-0.4cm}\left(p+\hdots+p^{(M-1)}+(N-M)p^M\right)(n+m)+m+np^M\\
  =\mathcal{O}\left((n+m)Np^M\right)\;.
\end{multline}
This formula aligns with the case of nominal MPC, where $p^M=1$ since $p=1$ and $M=0$. The number of equality constrains in \eqref{Equ:DRMPCProblem} is $(M+1)pn+n$, which can be reduced to $(p+1)n$ when using the sequential approach.

Clearly, the most significant factor for the computational complexity of DRMPC is the exponential term $p^M$ in \eqref{Equ:IneqConstraints}. It is due to the combination of all the corner cases in the tightening of the input and state constraints \eqref{Equ:InputTightening} and \eqref{Equ:StateTightening}. To put this into perspective, recall that the deadbeat horizon $M$ is selected as any number $M\geq n$. The number of states $n$ is typically small for Robust MPC applications (usually, $n\leq 4$). Furthermore, $p$ is the number of vertices of the disturbance set $\BD\subseteq\BR^n$, which may again be small, especially if $n$ is small and $\BD$ is not fully dimensional. For example, consider a system with $n=3$ states, where the disturbance is a hypercube affecting only 2 of the states. Then the DRMPC contains about $p^M=4^3=64$ times the number of constraints as the corresponding nominal MPC problem. In conclusion, DRMPC scales poorly with the dimension of the system state and the vertices of the disturbance set. It may be tractable for small systems. To circumvent the complexity problem, the deadbeat input sequences can be computed offline, as discussed next.

\vspace*{-0.4cm}
\subsection{Offline DRMPC}\label{Sec:OfflineDRMPC}

The case where the deadbeat input sequences are computed online (Online DRMPC) is the most general, but mainly of theoretical interest, due to its computational complexity. Presumably of greater practical interest is the case where they are pre-computed offline (Offline DRMPC). This significantly reduces the computational complexity and also allows for a less restrictive terminal condition. However, it may lead to sub-optimality of the solution to \eqref{Equ:DRMPCProblem}, and hence to a possible decrease in closed-loop performance. The offline pre-computations work as described below.

\begin{remark}[offline pre-computations]\label{Rem:Precomputation}
	(a) For each vertex $d^{(i)}$, the deadbeat control inputs $w_{0}^{(i)},\hdots,w_{M-1}^{(i)}$ and corresponding states $z_{0}^{(i)},\hdots,z_{M-1}^{(i)}$ are pre-computed and substituted as fixed values into the DRMPC problem \eqref{Equ:DRMPCProblem}. This eliminates them as decision variables, together with the contraints (\ref{Equ:DRMPCProblem}b,c). (b) To this end, a finite-horizon optimal control problem over the deadbeat horizon $M$ is solved for each vertex $d^{(i)}$. The vertex represents the initial condition, the nominal input and state constraints are imposed, and an analogous cost function to (\ref{Equ:DRMPCProblem}a) can be used.
\end{remark}

The pre-computed deadbeat control inputs and states are used to define the convex hulls
\begin{subequations}\label{Equ:OfflineDRMPCConvHull}\begin{align}
	\mathbb{W}_{k}&:=\conv{0,w^{(1)}_{k},w^{(2)}_{k},\dots,w^{(p)}_{k}}\;,\\
	\mathbb{Z}_{k}&:=\conv{0,z^{(1)}_{k},z^{(2)}_{k},\dots,z^{(p)}_{k}}
\end{align}\end{subequations}
for all $k=0,1,\dots,M-1$. Based on this, the tightened input constraints (\ref{Equ:DRMPCProblem}f,h) can be expressed as
\begin{subequations}\label{Equ:OfflineInputTightening}\begin{align}
  &u_{0|0}\in\BU\;,\\
  &u_{1|0}\in\BU\ominus\mathbb{W}_{0}\;,\\
  &\;\;\vdots\nonumber\\
  &u_{k|0}\in\BU\ominus\mathbb{W}_{0}\ominus\mathbb{W}_{1}\ominus\dots\ominus\mathbb{W}_{\max\{k-1,M-1\}}
\end{align}\end{subequations}
for all $k=0,1,\dots,N-1$. Analogously, the tightened state constraints (\ref{Equ:DRMPCProblem}g,i) can be expressed as
\begin{subequations}\label{Equ:OfflineInputTightening}\begin{align}
  &x_{1|0}\in\BX\ominus\mathbb{Z}_{0}\;,\\
  &\;\;\vdots\nonumber\\
  &x_{k|0}\in\BX\ominus\mathbb{Z}_{0}\ominus\mathbb{Z}_{1}\ominus\dots\ominus\mathbb{Z}_{\max\{k-1,M-1\}}
\end{align}\end{subequations}
for all $k=1,2,\dots,N$.

With these offline pre-computations, the complexity of Offline DRMPC will be similar to Nominal MPC. In particular, it will be independent of the dimensions $n$ and $p$, and of the choice for the deadbeat horizon $M$. All system-theoretic guarantees of Online DRMPC, which are discussed in Section \ref{Sec:Guarantees}, will be preserved.

\begin{remark}[alternative terminal conditions]\label{Rem:AltTerminalCondition}
  For Offline DRMPC, where the deadbeat input sequences are pre-computed offline, any terminal condition of a nominal MPC can be used, such as a positive invariant (PI) set.
\end{remark}

For example, if the stage cost function in (\ref{Equ:DRMPCProblem}a) is quadratic,
\begin{equation*}
	\ell\left(u_{k},x_{k}\right)=u_{k}\tp Ru_{k}+x_{k}\tp Qx_{k}
\end{equation*}
for some $R\succ 0$ and $Q\succeq 0$, the Linear Quadratic Regulator (LQR) $u_{k}=K_{\text{LQR}}x_{k}$ can be selected as the terminal controller. Then the terminal set for $x_{N|0}$ can be selected as any PI set for the LQR and
\begin{align*}
	u_{k}\in\BU\ominus\mathbb{W}_{0}\ominus\mathbb{W}_{1}\ominus\dots\ominus\mathbb{W}_{M-1}\;,\\
	x_{k}\in\BX\ominus\mathbb{Z}_{0}\ominus\mathbb{Z}_{1}\ominus\dots\ominus\mathbb{Z}_{M-1}\;.
\end{align*}
The terminal cost is then selected as $x_{N|0}\tp Px_{N|0}$, where $P\succeq 0$ is the solution to the Algebraic Riccati Equation of the LQR. 

Note that using an alternative terminal condition according to Remark \ref{Rem:AltTerminalCondition} will increase the region of attraction. However, a major point of this paper that no invariant set has to be calculated, even though in the case of DRMPC it may be a simpler (non-robust) PI set.

\section{System-theoretic Guarantees}\label{Sec:Guarantees}

It is shown in this section that DRMPC comes with all relevant system-theoretic guarantees, including recursive feasibility and stability. The line of argument is closely related to those of other RMPC approaches \cite{LimonEtAl:2009,RaMaDi:2018}.

\vspace*{-0.2cm}
\subsection{Recursive Feasibility}\label{Sec:Feasibility}

Under the DRMPC regime, in each time step $k=0,1,\hdots$ the first element of the optimal input sequence from \eqref{Equ:DRMPCProblem}, $u_{k}=u_{k|k}^{\star}$, is applied to the system.

\begin{theorem}[recursive feasibility]\label{The:RecFeasibility}
  If the DRMPC problem \eqref{Equ:DRMPCProblem} is feasible at $k=0$ for $\bar{x}_{0}$, it remains feasible for all future states $x_{1},x_{2},\hdots$ of system \label{Equ:DTSystem} under the DRMPC regime.
\end{theorem}

\emph{\textbf{Proof:}} Consider the solution to \eqref{Equ:DRMPCProblem} in step $k=0$. Let $u_{0|0}^{\star},\hdots,u_{N-1|0}^{\star}$ be the optimal inputs and $x_{0|0}^{\star},\hdots,x_{N|0}^{\star}$ be the optimal states. They satisfy the input constraints (\ref{Equ:DRMPCProblem}h) and state constraints (\ref{Equ:DRMPCProblem}i), respectively. By the terminal constraint (\ref{Equ:DRMPCProblem}e), $x_{N|0}^{\star}=0$. Moreover, let $w_{0}^{(i)\star},\hdots,w_{M-1}^{(i)\star}$ be the optimal deadbeat inputs and $z_{0}^{(i)\star},\hdots,z_{M}^{(i)\star}$ be the corresponding states, with $z_{M}^{(i)\star}=0$, for any $i=1,\hdots,p$.

For recursive feasibility, it must be shown that in time step $k=1$ there exists a candidate input sequence $u_{1|1},\hdots,u_{N|1}$ and a corresponding state sequence $x_{1|1},\hdots,x_{N+1|1}$ that satisfy the input, state, and terminal constraints, for any admissible disturbance $d_{0}$. For the sake of this proof, it is assumed that $w_{0}^{(i)\star},\hdots,w_{M-1}^{(i)\star}$ and $z_{0}^{(i)\star},\hdots,z_{M}^{(i)\star}$ remain unchanged from step 0 to step 1, which is a feasible choice.

By the presence of a disturbance,
\begin{equation*}
  x_{1|1}=x_{1|0}^{\star}+d_{0}\;,\quad\text{for some}\;d_{0}\in\BD\;.
\end{equation*}
Thus $d_{0}$ is known in step 1, and according to \eqref{Equ:DisConvHull} it can be expressed as\vspace*{-0.2cm}
\begin{equation*}
  d_{0}=\sum_{i=1}^{p}\lambda_{0}^{(i)}d^{(i)}\;,
\end{equation*}
for some $\lambda_{0}^{(i)}\geq 0$ with $\sum_{i=1}^{p}\lambda_{0}^{(i)}=1$. Consider the candidate input sequence\vspace*{-0.2cm}
\begin{align*}
  u_{1|1}&=u_{1|0}^{\star}+\sum_{i=1}^{p}\lambda_{0}^{(i)}w_{0}^{(i)\star}\;,\\
  &\;\;\vdots\\
  u_{M|1}&=u_{M|0}^{\star}+\sum_{i=1}^{p}\lambda_{0}^{(i)}w_{M-1}^{(i)\star}\;,
\end{align*}

\begin{align*}  
  u_{M+1|1}&=u_{M+1|0}^{\star}\;,\\
  &\;\;\vdots\\
  u_{N-1|1}&=u_{N-1|0}^{\star}\;,\\
  u_{N|1}&=0\;.\\
\end{align*}\vspace*{-0.8cm}

\noindent By construction, this input sequence satisfies all the input constraints (\ref{Equ:DRMPCProblem}h). Moreover, the corresponding state sequence\vspace*{-0.1cm}
\begin{align*}
  x_{1|1}&=x_{1|0}^{\star}+\sum_{i=1}^{p}\lambda_{0}^{(i)}z_{0}^{(i)\star}\;,\\
  &\;\;\vdots\\
  x_{M+1|1}&=x_{M+1|0}^{\star}+\sum_{i=1}^{p}\lambda_{0}^{(i)}z_{M}^{(i)\star}\;,\\
  x_{M+2|1}&=x_{M+2|0}^{\star}\;,\\
  &\;\;\vdots\\
  x_{N|1}&=x_{N|0}^{\star}\;,\\
  x_{N+1|1}&=0\;.\\
\end{align*}\vspace*{-0.8cm}

\noindent satisfies all the state constraints (\ref{Equ:DRMPCProblem}i) and the terminal constraint. The last step follows from $x_{N|0}^{\star}=0$ and the choice of $u_{N|1}=0$.\hfill$\square$

\vspace*{-0.2cm}
\subsection{Input-to-State Stability}\label{Sec:Stability}

Let $\CX_{N}\subseteq\BX$ be the \emph{feasible set}, i.e., the set of initial conditions $x_{0}$ for which problem \eqref{Equ:DRMPCProblem} is feasible. The \emph{state feedback control law} $\kappa_{N}:\CX_{N}\to\BU$ of the DRMPC regime is defined as $\kappa_{N}(x_{k})=u_{k|k}^{\star}$. From Theorem \ref{The:RecFeasibility}, this map is well-defined.

\begin{definition}[robust positive invariant set]
  Consider system \eqref{Equ:DTSystem} with any state feedback control law $\kappa:\BX\to\BU$. A set $\CP\subseteq\BX$ is called a \textbf{robust positive invariant (RPI) set} if for all $x\in\CP$ it holds that
  \begin{equation*}
    \kappa(x)\in\BU\quad\text{and}\quad Ax+B\kappa(x)\oplus\BD\subseteq\CP\;.
  \end{equation*}
\end{definition}

\begin{corollary}[invariance of the feasible set]\label{The:InvFeasSet}
  The feasible set $\CX_{N}$ is a RPI set for system \eqref{Equ:DTSystem} under the DRMPC regime.
\end{corollary}

\emph{\textbf{Proof:}} This is immediate from Theorem \ref{The:RecFeasibility}.\hfill$\square$

Computing the feasible set of the DRMPC problem \eqref{Equ:DRMPCProblem} hence leads to an RPI set. Unlike most other RPI sets, this one is based on a \emph{nonlinear} disturbance feedback policy, following Remark \ref{Rem:NonlinearDF}.

On top of this invariance property, it can be shown that the DRMPC stabilizes system \eqref{Equ:DTSystem}. Due to the presence of disturbances, the appropriate stability concept to be used is \emph{input-to-state stability (ISS)} \cite{JiangWang:2001,LimonEtAl:2009,RaMaDi:2018}. It requires the following standard assumption on the stage cost function.

\begin{assumption}[stage cost]\label{Ass:StageCost}
  For the stage cost function $\ell:\BX\times\BU\rightarrow\BR_{0+}$ there exist two $K_{\infty}$-functions $\varepsilon_{1}$ and $\epsilon_{2}$ with
  \begin{equation}\label{Equ:LyaFunction1}
    \varepsilon_{1}(\|x\|)\leq \ell\left(x,u\right)\leq \varepsilon_{2}(\|x\|)\quad\fa x\in\BX\;,\;u\in\BU\;.
  \end{equation}
\end{assumption}

\begin{definition}[input-to-state stability]\label{Def:Stability}
    System \eqref{Equ:DTSystem} under some state feedback law $\kappa:\CP\to\BU$ is \textbf{input-to-state stable (ISS)} on some RPI set $\CP\subseteq\BX$ if there exist a $\text{KL}$-function $\beta$ and a $\text{K}$-function $\gamma$ such that
    \begin{equation*}
        \|x_{k}\|\leq\beta\bigl(\|x_{0}\|,k\bigr)+\gamma\bigl(\max_{j<k}\|d_{j}\|\bigr) \quad\fa k\in\BZ_{0+}\;,
    \end{equation*}
    for any $x_{0}\in\CP$ and any disturbances $d_{0},d_{1},\hdots\in\BD$.
\end{definition}

Definition \ref{Def:Stability} reduces to \emph{asymptotic stability} if all disturbances vanish beyond some step $k$, i.e., $d_{j}=0$ for all $j\geq k$. If the disturbances are non-zero, the states $x_{k}$ are bounded by norm balls, whose size increases with the maximum value of $\|d_{j}\|$ observed in the past $j<k$.

\begin{definition}[ISS Lyapunov function]\label{Def:LyaFunction}
  $V:\CP\rightarrow\BR_{0+}$ is called an \textbf{ISS Lyapunov function} for system \eqref{Equ:DTSystem} under the state feedback law $\kappa:\CP\to\BU$, if the following conditions are satisfied:\\
    (a) There exist two $K_{\infty}$-functions $\alpha_{1}$ and $\alpha_{2}$ such that
    \begin{equation}\label{Equ:LyaFunction1}
        \alpha_{1}(\|x\|)\leq V\left(x\right)\leq \alpha_{2}(\|x\|)\quad\fa x\in\CP\;.
    \end{equation}
    (b) There exists a $K_{\infty}$-function $\alpha_{3}$ and $\text{K}$-function $\sigma$ with
    \begin{equation}\label{Equ:LyaFunction2}
        V\left(x_{k+1}\right)\leq V\left(x_{k}\right)-\alpha_{3}(\|x_{k}\|)+\sigma(\|d_{k}\|)
    \end{equation}
    for all $k\in\BZ_{0+}$, any $x_{0}\in\CP$, and any $d_{k}\in\BD$.
\end{definition}

Definition \ref{Def:LyaFunction} reduces to a \emph{Lyapunov function} if the term $\sigma(\|d_{k}\|)$ is removed in \eqref{Equ:LyaFunction2}.

\begin{lemma}[Lyapunov stability]\label{The:LyapunovStability}
  If there exists a ISS Lyapunov function $V:\CP\rightarrow\BR_{0+}$ for system \eqref{Equ:DTSystem} under some feedback law $\kappa:\BX\to\BU$ and on some RPI set $\CP\subseteq\BX$, then the closed-loop system is ISS.
\end{lemma}

\emph{\textbf{Proof:}} See \cite{JiangWang:2001}.

\begin{theorem}[input-to-state stability]\label{The:Stability}
  The closed-loop system \eqref{Equ:DTSystem} under the DRMPC regime $\kappa_{N}:\CX_{N}\to\BU$ is ISS on $\CX_{N}$.
\end{theorem}

\emph{\textbf{Proof:}} By Lemma \ref{The:LyapunovStability}, it suffices to show that the optimal value function of the DRMPC problem $V^{\star}_{N}:\CX_{N}\to\BR_{0+}$ is an ISS Lyapunov function in the sense of Definition \ref{Def:LyaFunction}.

Definition \ref{Def:LyaFunction}(a) is satisfied by virtue of Assumption \ref{Ass:StageCost}. Definition \ref{Def:LyaFunction}(b) is verified as follows. 

First, for the case of $d_{k}=0$, the shifted optimal input sequence from the previous step
\begin{align*}
  u_{k+1|k+1}&=u^{\star}_{k+1|k}\;,\\
  &\dots\;,\\
  u_{k+N|k+1}&=u^{\star}_{k+N|k}\;,\\
  u_{k+N+1|k+1}&=0
\end{align*}
is feasible, and so is the corresponding state sequence. Note that this is due to the terminal constraint (\ref{Equ:DRMPCProblem}e) as well as the constraints (\ref{Equ:DRMPCProblem}f,g) and the zeros in (\ref{Equ:OfflineDRMPCConvHull}a,b), which are for this reason essential. By the standard argument in MPC stability proofs, the optimal value function hence decreases by more than the first stage cost,
\begin{equation*}
  V^{\star}_{N}\left(x_{k+1}\right)\leq V^{\star}_{N}\left(x_{k}\right)-\ell\left(x_{k},u_{k}\right)\;.
\end{equation*}
Together with Assumption \ref{Ass:StageCost} this verifies the term $-\alpha_{3}(\|x_{k}\|)$ in \eqref{Equ:LyaFunction2}, for the case of $d_{k}=0$.

Second, for the case of $d_{k}\neq 0$, observe that the optimal value function $V^{\star}_{N}(x_{k+1|k}^{\star}+d_{k})$ is continuous in $d_{k}$, for all $d_{k}\in\mathbb{D}$. Together with the compactness of $\mathbb{D}$, this verifies the term $+\sigma(\|d_{k}\|)$ in \eqref{Equ:LyaFunction2}, for the case of $d_{k}\neq 0$.\hfill$\square$

\section{Numerical Example}\label{Sec:Example}

For the purposes of illustration and evaluation, a small numerical example is considered.

\vspace*{-0.4cm}
\subsection{Example System}

Consider a discrete-time double integrator with additive disturbances,
\begin{equation}
  x_{k+1}=
  \begin{bmatrix}
    1 & t_{\mathrm{S}} \\
    0 & 1
  \end{bmatrix}x_{k}+
  \begin{bmatrix}
    t_{\mathrm{S}}^2/2 \\
    t_{\mathrm{S}}
  \end{bmatrix}u_{k}+d_{k}\;.
\end{equation}
Note that the sampling time $t_{\mathrm{S}}$ is a parameter of the system. The state and input constraints are given by
\begin{subequations}\begin{align}
  \BX&=\left\{x\in\BR^2\;:\;
    \begin{bmatrix}
      -5\\
      -2
    \end{bmatrix}\leq x\leq
    \begin{bmatrix}
      +5\\
      +2
    \end{bmatrix}\right\}\;,\\
  \BU&=\left\{u\in\BR\;:\;
    -2\leq u\leq +2\right\}\;. 
\end{align}\end{subequations}
The disturbances are contained in the disturbance set
\begin{equation}
  \BD=\left\{d\in\BR^2\;:\;
    \begin{bmatrix}
      -0.15\\
      -0.15
    \end{bmatrix}\leq d\leq
    \begin{bmatrix}
      +0.15\\
      +0.15
    \end{bmatrix}\right\}\;.
\end{equation}
It has the following $p=4$ vertices: $d^{(i)}=[\pm0.15,\pm0.15]\tp$. The stage cost function is assumed as a quadratic function
\begin{subequations}\label{Equ:ExCostFunction}\begin{align}
  &\ell\left(x,u\right)=x\tp Qx+u\tp Ru\;,\\
  \hspace*{-1cm}\text{with}\quad
  &Q=\begin{bmatrix}
   1 & 0\\
   0 & 1
  \end{bmatrix}\;,\quad
  R=10\;.
\end{align}\end{subequations}

\vspace*{-0.4cm}
\subsection{Controller Evaluation}

All of the RMPC considered below come with the guarantee of (robust) recursive feasibility and stability. Hence the next most relevant performance criterion for RMPC is the \textbf{region of attraction (RoA)}. The RoA is the set of initial states that can hence be stabilized, i.e., for which the respective MPC problem is feasible.

\subsubsection*{Comparison of Online and Offline DRMPC}

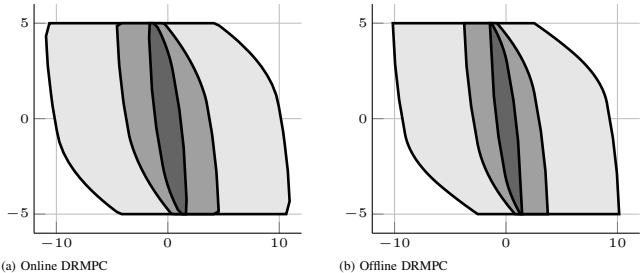
\begin{figure}[H]
	\hspace*{-0.3cm}
    \tiny
	\begin{tabular}{p{0.46\columnwidth}p{0.46\columnwidth}}

\begin{minipage}[t]{0.4\columnwidth}

\begin{tikzpicture}

\begin{axis}[%
width=1.400in,
height=1.200in,
at={(0.0in,0.0in)},
scale only axis,
xmin=-12,
xmax=12,
ymin=-6,
ymax=6,
label style={font=\tiny},
tick label style={font=\tiny},
axis background/.style={fill=white},
axis x line*=bottom,
axis y line*=left,
xmajorgrids,
ymajorgrids
]


\addplot[area legend, line width=1.0pt, draw=black, fill=LiteGray, forget plot]
table[row sep=crcr] {
x	y\\
0.0000 -5.0000\\
9.0950 -5.0000\\
10.6252 -4.9998\\
10.9269 -4.3263\\
10.7705 -2.7654\\
10.4891 -1.3251\\
10.3382 -0.6504\\
9.9608 0.6267\\
9.7282 1.2290\\
9.3593 1.7854\\
8.9027 2.2858\\
8.3970 2.7283\\
7.8733 3.1173\\
7.3518 3.4595\\
6.8361 3.7582\\
6.3545 4.0327\\
5.8841 4.2751\\
5.4334 4.4949\\
5.0010 4.6963\\
4.5843 4.8818\\
4.1360 4.9995\\
3.6327 5.0000\\
0.0000 5.0000\\
-9.0950 5.0000\\
-10.6252 4.9998\\
-10.9275 4.3265\\
-10.8509 3.5257\\
-10.7705 2.7654\\
-10.4891 1.3251\\
-10.3382 0.6504\\
-9.9614 -0.6267\\
-9.7282 -1.2290\\
-9.3593 -1.7854\\
-8.9021 -2.2857\\
-8.3970 -2.7283\\
-7.8733 -3.1173\\
-7.3518 -3.4595\\
-6.8449 -3.7630\\
-6.3540 -4.0323\\
-5.8841 -4.2751\\
-5.4334 -4.4949\\
-5.0010 -4.6963\\
-4.5843 -4.8818\\
-4.1360 -4.9995\\
-3.6327 -5.0000\\
}--cycle;


\addplot[area legend, line width=1.0pt, draw=black, fill=MedGray, forget plot]
table[row sep=crcr] {
x	y\\
0.0000 -4.7596\\
0.3136 -4.9845\\
0.6316 -4.9999\\
3.6327 -5.0000\\
4.1360 -4.9995\\
4.5628 -4.8589\\
4.4513 -3.6824\\
4.4044 -3.2000\\
4.3550 -2.7637\\
4.2746 -2.3500\\
4.2019 -1.9773\\
4.0743 -1.3238\\
4.0172 -1.0314\\
3.9535 -0.7542\\
3.8757 -0.4896\\
3.5326 0.6739\\
3.4644 0.8895\\
3.3784 1.0977\\
3.2794 1.2984\\
3.1704 1.4919\\
3.0577 1.6810\\
2.9393 1.8653\\
2.8143 2.0447\\
2.6871 2.2230\\
2.5505 2.3951\\
2.4123 2.5689\\
2.2647 2.7376\\
2.1124 2.9074\\
1.9536 3.0784\\
1.6108 3.4231\\
1.4235 3.5953\\
1.2266 3.7751\\
1.0161 3.9574\\
0.7908 4.1457\\
0.5485 4.3420\\
0.2859 4.5449\\
0.0000 4.7596\\
-0.3136 4.9845\\
-0.6316 4.9999\\
-3.6327 5.0000\\
-4.1360 4.9995\\
-4.5628 4.8589\\
-4.4513 3.6824\\
-4.4055 3.2007\\
-4.3550 2.7637\\
-4.2746 2.3500\\
-4.2019 1.9773\\
-4.0743 1.3238\\
-4.0172 1.0314\\
-3.9535 0.7542\\
-3.8757 0.4896\\
-3.5326 -0.6739\\
-3.4644 -0.8895\\
-3.3784 -1.0977\\
-3.2788 -1.2982\\
-3.1704 -1.4919\\
-3.0572 -1.6807\\
-2.9393 -1.8653\\
-2.8143 -2.0447\\
-2.6871 -2.2230\\
-2.5505 -2.3951\\
-2.4123 -2.5689\\
-2.2647 -2.7376\\
-2.1124 -2.9074\\
-1.9536 -3.0784\\
-1.6108 -3.4231\\
-1.4235 -3.5953\\
-1.2266 -3.7751\\
-1.0161 -3.9574\\
-0.7908 -4.1457\\
-0.5485 -4.3420\\
-0.2859 -4.5449\\
}--cycle;


\addplot[area legend, line width=1.0pt, draw=black, fill=DarkGray, forget plot]
table[row sep=crcr] {
x	y\\
0.0000 -3.5070\\
0.2460 -3.9107\\
0.5476 -4.3345\\
0.9118 -4.7798\\
1.2837 -4.9995\\
1.6244 -4.9994\\
1.6803 -4.2441\\
1.5440 -2.8085\\
1.4693 -2.0223\\
1.4209 -1.5131\\
1.3787 -1.2947\\
1.2615 -0.6935\\
1.0725 0.2754\\
1.0315 0.4854\\
0.5495 2.1402\\
0.4637 2.4307\\
0.3493 2.7646\\
0.1966 3.1246\\
0.0000 3.5070\\
-0.2460 3.9107\\
-0.5476 4.3345\\
-0.9118 4.7798\\
-1.2837 4.9995\\
-1.6244 4.9994\\
-1.6803 4.2441\\
-1.5440 2.8085\\
-1.4693 2.0223\\
-1.4209 1.5131\\
-1.3787 1.2947\\
-1.2615 0.6935\\
-1.0591 -0.3441\\
-1.0315 -0.4854\\
-0.5495 -2.1402\\
-0.4637 -2.4307\\
-0.3493 -2.7646\\
-0.1966 -3.1246\\
}--cycle;

\end{axis}

\end{tikzpicture}

\end{minipage}
    &

\begin{minipage}[t]{0.4\columnwidth}

\begin{tikzpicture}

\begin{axis}[%
width=1.400in,
height=1.200in,
at={(0.0in,0.0in)},
scale only axis,
xmin=-12,
xmax=12,
ymin=-6,
ymax=6,
label style={font=\tiny},
tick label style={font=\tiny},
axis background/.style={fill=white},
axis x line*=bottom,
axis y line*=left,
xmajorgrids,
ymajorgrids
]


\addplot[area legend, line width=1.0pt, draw=black, fill=LiteGray, forget plot]
table[row sep=crcr] {
x	y\\
10.1550 -5.0000\\
9.9507 -2.8500\\
9.5900 -1.0000\\
9.0443 0.8495\\
8.9816 1.0084\\
8.9029 1.1673\\
8.8084 1.3262\\
8.6979 1.4851\\
8.5716 1.6440\\
8.4294 1.8029\\
8.2713 1.9618\\
8.0973 2.1207\\
7.9074 2.2796\\
7.7016 2.4385\\
7.4800 2.5974\\
7.2424 2.7563\\
6.9890 2.9152\\
6.7196 3.0741\\
6.4344 3.2330\\
6.1333 3.3919\\
5.7575 3.5803\\
5.0917 3.8981\\
4.3941 4.2159\\
3.6647 4.5337\\
2.9036 4.8515\\
2.5331 5.0000\\
-10.1550 5.0000\\
-9.9507 2.8500\\
-9.5900 1.0000\\
-9.0443 -0.8495\\
-8.9816 -1.0084\\
-8.9029 -1.1673\\
-8.8084 -1.3262\\
-8.6979 -1.4851\\
-8.5716 -1.6440\\
-8.4294 -1.8029\\
-8.2713 -1.9618\\
-8.0973 -2.1207\\
-7.9074 -2.2796\\
-7.7016 -2.4385\\
-7.4800 -2.5974\\
-7.2424 -2.7563\\
-6.9890 -2.9152\\
-6.7196 -3.0741\\
-6.4344 -3.2330\\
-6.1333 -3.3919\\
-5.7575 -3.5803\\
-5.0917 -3.8981\\
-4.3941 -4.2159\\
-3.6647 -4.5337\\
-2.9036 -4.8515\\
-2.5331 -5.0000\\
}--cycle;


\addplot[area legend, line width=1.0pt, draw=black, fill=MedGray, forget plot]
table[row sep=crcr] {
x	y\\
0.1489 -4.5337\\
0.5922 -4.8515\\
0.8142 -5.0000\\
3.7526 -5.0000\\
3.5396 -2.7583\\
3.1290 -0.6527\\
2.5365 1.3557\\
2.4110 1.6735\\
2.2537 1.9913\\
2.0646 2.3091\\
1.8437 2.6269\\
1.5911 2.9447\\
1.3066 3.2625\\
0.9904 3.5803\\
0.6424 3.8981\\
0.2627 4.2159\\
-0.1489 4.5337\\
-0.5922 4.8515\\
-0.8142 5.0000\\
-3.7526 5.0000\\
-3.5396 2.7583\\
-3.1290 0.6527\\
-2.5365 -1.3557\\
-2.4110 -1.6735\\
-2.2537 -1.9913\\
-2.0646 -2.3091\\
-1.8437 -2.6269\\
-1.5911 -2.9447\\
-1.3066 -3.2625\\
-0.9904 -3.5803\\
-0.6424 -3.8981\\
-0.2627 -4.2159\\
}--cycle;


\addplot[area legend, line width=1.0pt, draw=black, fill=DarkGray, forget plot]
table[row sep=crcr] {
x	y\\
0.0440 -3.2625\\
0.2013 -3.5803\\
0.3904 -3.8981\\
0.6113 -4.2159\\
0.8639 -4.5337\\
1.1484 -4.8515\\
1.2961 -5.0000\\
1.4485 -5.0000\\
1.0846 -1.1693\\
0.6740 0.9363\\
0.0815 2.9447\\
-0.0440 3.2625\\
-0.2013 3.5803\\
-0.3904 3.8981\\
-0.6113 4.2159\\
-0.8639 4.5337\\
-1.1484 4.8515\\
-1.2961 5.0000\\
-1.4485 5.0000\\
-1.0846 1.1693\\
-0.6740 -0.9363\\
-0.0815 -2.9447\\
}--cycle;

\end{axis}

\end{tikzpicture}

\end{minipage}\\
    (a) Online DRMPC & (b) Offline DRMPC 
    \end{tabular}
    \normalsize
	\caption{RoA for $t_{\mathrm{S}}=0.1$ and $N=10$ (dark gray), $N=20$ (medium gray), $N=40$ (light gray). The deadbeat horizon is $M=3$ and the terminal set is the origin in both cases.\label{Fig:RoADRMPConoff}}
\end{figure}

The first part of the numerical example compares \textbf{Online DRMPC} against \textbf{Offline DRMPC}. Consider Problem \eqref{Equ:DRMPCProblem} and Remark \ref{Rem:Precomputation} for the details. Figure \ref{Fig:RoADRMPConoff} shows a direct comparison of the respective RoA for different horizon lengths. 

As expected, the Online DRMPC always results in a larger RoA than Offline DRMPC. A quantitative analysis of the excess volume of the RoA of Online DRMPC over that of Offline DRMPC is shown in Table \ref{Tab:RoADRMPConoff}, for different values of the sampling time $t_{\mathrm{S}}$ and horizon length $N$. Observe that this excess volume of the RoA converges to zero for large horizon lengths (and large values of the parameter $t_{\mathrm{S}}$).

\renewcommand{\arraystretch}{1.4}
\begin{table}[H]
	\centering
	\begin{tabular}{|l|l||D{.}{.}{3.2}|D{.}{.}{3.2}|D{.}{.}{3.2}|D{.}{.}{3.2}|D{.}{.}{3.2}|}
		\hline
		\multicolumn{2}{|c||}{$\,$} & \multicolumn{5}{c|}{\textbf{horizon length $N$}}\\\cline{3-7}
		\multicolumn{2}{|c||}{$\,$} & \multicolumn{1}{c|}{$10$} & \multicolumn{1}{c|}{$20$} & \multicolumn{1}{c|}{$30$} & \multicolumn{1}{c|}{$40$} & \multicolumn{1}{c|}{$50$}\\\hline\hline
		\multirow{4}{0.3cm}{\rotatebox{90}{\textbf{parameter $t_{\mathrm{S}}$}}}
		 	& $0.10$ & +37.9\% & +27.6\% & +16.7\% &+11.7\% & +8.2\%\\\cline{2-7}
			& $0.15$ &  +6.6\% &  +2.5\% &  +1.3\% & +0.0\% & +0.0\%\\\cline{2-7}
			& $0.25$ &  +2.8\% &  +1.0\% &  +0.0\% & +0.0\% & +0.0\%\\\cline{2-7}
			& $0.40$ &  +2.1\% &  +0.0\% &  +0.0\% & +0.0\% & +0.0\%\\\hline
	\end{tabular}
	\caption{Excess volume of RoA of Online DRMPC over Offline DRMPC\label{Tab:RoADRMPConoff}}
\end{table}
\renewcommand{\arraystretch}{1.0}

Offline DRMPC, however, is computationally simpler than Online DRMPC. Given a fixed amount of computation time, a longer horizon can hence be chosen for Offline DRMPC, which again increases its RoA. There is no quantitative analysis of computation times here, because it is highly dependent (i) on individual problem data and dimensions and (ii) on the selected optimization solver and hardware setup.

\subsubsection*{Comparison with Other RMPC Approaches}

The second part of the numerical example compares Offline DRMPC with competitive RMPC approaches, namely Tube-based RMPC (`Tube-RMPC') and RMPC with Affine Disturbance Feedback (`RMPC-ADF').

\textbf{DRMPC:} The deadbeat input sequences of DRMPC are computed offline, according to Remark \ref{Rem:Precomputation}. The terminal constraint is the maximal PI set. The deadbeat horizon is selected as $M=N$ in all cases.

\textbf{Tube-RMPC:} The Tube-RMPC setup follows the formulation in \cite{Mayne:2005}. The system is pre-stabilized using the LQR feedback gain $K_{\text{LQR}}$ corresponding to \eqref{Equ:ExCostFunction}. The tubes are calculated as the minimal RPI set as proposed in \cite{RacEtAl:2005}.

\textbf{RMPC-ADF:} The disturbance feedback gain is fixed offline to $K_{\text{ADF}}=K_{\text{LQR}}$. This reduces the computational complexity of Affine Disturbance Feedback \cite{Goulart:2006}, such that it represents a fair comparison to Offline DRMPC.

\begin{figure}[H]
    \tiny
	\begin{tabular}{p{0.46\columnwidth}p{0.46\columnwidth}}

\begin{minipage}[t]{0.4\columnwidth}

\begin{tikzpicture}

\begin{axis}[%
width=1.400in,
height=1.200in,
at={(0.0in,0.0in)},
scale only axis,
xmin=-20,
xmax=20,
ymin=-6,
ymax=6,
label style={font=\tiny},
tick label style={font=\tiny},
axis background/.style={fill=white},
axis x line*=bottom,
axis y line*=left,
xmajorgrids,
ymajorgrids
]


\addplot[area legend, line width=1.0pt, draw=black, fill=LiteGray, forget plot]
table[row sep=crcr] {
x	y\\
17.3436 -5.0000\\
17.4693 -4.9125\\
17.5177 -4.7957\\
17.5472 -4.6025\\
17.5472 -4.3025\\
17.3924 -2.6730\\
17.0746 -1.0435\\
16.5939 0.5860\\
16.3050 1.3175\\
16.2140 1.4096\\
16.1258 1.4987\\
16.0311 1.5944\\
15.9268 1.6996\\
15.8376 1.7894\\
15.7275 1.8999\\
15.6500 1.9774\\
15.5770 2.0502\\
15.4961 2.1307\\
15.4244 2.2019\\
15.3268 2.2983\\
15.2186 2.4045\\
15.0988 2.5215\\
14.9662 2.6498\\
14.8196 2.7904\\
14.6577 2.9436\\
14.4793 3.1100\\
14.2830 3.2895\\
14.0675 3.4819\\
13.8317 3.6860\\
13.5745 3.8997\\
13.2954 4.1195\\
12.9943 4.3400\\
12.6718 4.5529\\
12.3299 4.7464\\
11.9721 4.9031\\
11.6044 4.9990\\
-17.3436 5.0000\\
-17.4693 4.9125\\
-17.5177 4.7957\\
-17.5472 4.6025\\
-17.5472 4.3025\\
-17.3924 2.6730\\
-17.0746 1.0435\\
-16.5939 -0.5860\\
-16.3050 -1.3175\\
-16.2139 -1.4096\\
-16.1258 -1.4987\\
-16.0311 -1.5944\\
-15.9268 -1.6996\\
-15.8376 -1.7894\\
-15.7275 -1.8999\\
-15.6500 -1.9774\\
-15.5770 -2.0502\\
-15.4961 -2.1307\\
-15.4244 -2.2019\\
-15.3268 -2.2983\\
-15.2186 -2.4045\\
-15.0988 -2.5215\\
-14.9662 -2.6498\\
-14.8196 -2.7904\\
-14.6577 -2.9436\\
-14.4793 -3.1100\\
-14.2830 -3.2895\\
-14.0675 -3.4819\\
-13.8317 -3.6860\\
-13.5745 -3.8997\\
-13.2954 -4.1195\\
-12.9943 -4.3400\\
-12.6718 -4.5529\\
-12.3299 -4.7464\\
-11.9721 -4.9031\\
-11.6044 -4.9990\\
}--cycle;


\addplot[area legend, line width=1.0pt, draw=black, fill=MedGray, forget plot]
table[row sep=crcr] {
x	y\\
11.7236 -5.0000\\
11.8493 -4.9125\\
11.8977 -4.7957\\
11.9271 -4.6025\\
11.9271 -4.3025\\
11.7723 -2.6730\\
11.4546 -1.0435\\
10.9739 0.5860\\
10.6849 1.3175\\
10.5939 1.4096\\
10.5058 1.4987\\
10.4111 1.5944\\
10.3068 1.6996\\
10.2176 1.7894\\
10.1074 1.8999\\
10.0299 1.9774\\
9.9570 2.0502\\
9.8761 2.1307\\
9.8043 2.2019\\
9.7067 2.2983\\
9.5985 2.4045\\
9.4787 2.5215\\
9.3461 2.6498\\
9.1995 2.7904\\
9.0377 2.9436\\
8.8593 3.1100\\
8.6629 3.2895\\
8.4475 3.4819\\
8.2116 3.6860\\
7.9545 3.8997\\
7.6754 4.1195\\
7.3742 4.3400\\
7.0518 4.5529\\
6.7098 4.7464\\
6.3520 4.9031\\
5.9844 4.9990\\
-11.7236 5.0000\\
-11.8493 4.9125\\
-11.8977 4.7957\\
-11.9271 4.6025\\
-11.9271 4.3025\\
-11.7723 2.6730\\
-11.4546 1.0435\\
-10.9739 -0.5860\\
-10.6849 -1.3175\\
-10.5939 -1.4096\\
-10.5058 -1.4987\\
-10.4111 -1.5944\\
-10.3068 -1.6996\\
-10.2176 -1.7894\\
-10.1074 -1.8999\\
-10.0299 -1.9774\\
-9.9570 -2.0502\\
-9.8761 -2.1307\\
-9.8043 -2.2019\\
-9.7067 -2.2983\\
-9.5985 -2.4045\\
-9.4787 -2.5215\\
-9.3461 -2.6498\\
-9.1995 -2.7904\\
-9.0377 -2.9436\\
-8.8593 -3.1100\\
-8.6629 -3.2895\\
-8.4475 -3.4819\\
-8.2116 -3.6860\\
-7.9545 -3.8997\\
-7.6754 -4.1195\\
-7.3742 -4.3400\\
-7.0518 -4.5529\\
-6.7098 -4.7464\\
-6.3520 -4.9031\\
-5.9844 -4.9990\\
}--cycle;


\addplot[area legend, line width=1.0pt, draw=black, fill=DarkGray, forget plot]
table[row sep=crcr] {
x	y\\
8.9135 -5.0000\\
9.0392 -4.9125\\
9.0876 -4.7957\\
9.1171 -4.6025\\
9.1171 -4.3025\\
8.9623 -2.6730\\
8.6445 -1.0435\\
8.1638 0.5860\\
7.8749 1.3175\\
7.7839 1.4096\\
7.6958 1.4987\\
7.6011 1.5944\\
7.4968 1.6996\\
7.4075 1.7894\\
7.2974 1.8999\\
7.2199 1.9774\\
7.1470 2.0502\\
7.0660 2.1307\\
6.9943 2.2019\\
6.8967 2.2983\\
6.7885 2.4045\\
6.6687 2.5215\\
6.5361 2.6498\\
6.3895 2.7904\\
6.2276 2.9436\\
6.0492 3.1100\\
5.8529 3.2895\\
5.6374 3.4819\\
5.4016 3.6860\\
5.1445 3.8997\\
4.8654 4.1195\\
4.5642 4.3400\\
4.2417 4.5529\\
3.8998 4.7464\\
3.5420 4.9031\\
3.1743 4.9990\\
-8.9135 5.0000\\
-9.0392 4.9125\\
-9.0876 4.7957\\
-9.1171 4.6025\\
-9.1171 4.3025\\
-8.9623 2.6730\\
-8.6445 1.0435\\
-8.1638 -0.5860\\
-7.8749 -1.3175\\
-7.7839 -1.4096\\
-7.6958 -1.4987\\
-7.6011 -1.5944\\
-7.4968 -1.6996\\
-7.4075 -1.7894\\
-7.2974 -1.8999\\
-7.2199 -1.9774\\
-7.1470 -2.0502\\
-7.0660 -2.1307\\
-6.9943 -2.2019\\
-6.8967 -2.2983\\
-6.7885 -2.4045\\
-6.6687 -2.5215\\
-6.5361 -2.6498\\
-6.3895 -2.7904\\
-6.2276 -2.9436\\
-6.0492 -3.1100\\
-5.8529 -3.2895\\
-5.6374 -3.4819\\
-5.4016 -3.6860\\
-5.1445 -3.8997\\
-4.8654 -4.1195\\
-4.5642 -4.3400\\
-4.2417 -4.5529\\
-3.8998 -4.7464\\
-3.5420 -4.9031\\
-3.1743 -4.9990\\
}--cycle;

\end{axis}

\end{tikzpicture}

\end{minipage}
    &

\begin{minipage}[t]{0.4\columnwidth}

\begin{tikzpicture}

\begin{axis}[%
width=1.400in,
height=1.200in,
at={(0.0in,0.0in)},
scale only axis,
xmin=-20,
xmax=20,
ymin=-6,
ymax=6,
label style={font=\tiny},
tick label style={font=\tiny},
axis background/.style={fill=white},
axis x line*=bottom,
axis y line*=left,
xmajorgrids,
ymajorgrids
]


\addplot[area legend, line width=1.0pt, draw=black, fill=LiteGray, forget plot]
table[row sep=crcr] {
x	y\\
18.6999 -5.0000\\
18.4956 -2.8500\\
18.0972 -0.8068\\
17.5170 1.1599\\
16.7615 3.0726\\
15.8334 4.9474\\
-18.6999 5.0000\\
-18.4956 2.8500\\
-18.0972 0.8068\\
-17.5170 -1.1599\\
-16.7615 -3.0726\\
-15.8334 -4.9474\\
}--cycle;


\addplot[area legend, line width=1.0pt, draw=black, fill=MedGray, forget plot]
table[row sep=crcr] {
x	y\\
12.6887 -5.0000\\
12.4845 -2.8500\\
12.0860 -0.8068\\
11.5058 1.1599\\
10.7503 3.0726\\
9.8223 4.9474\\
-12.6887 5.0000\\
-12.4845 2.8500\\
-12.0860 0.8068\\
-11.5058 -1.1599\\
-10.7503 -3.0726\\
-9.8223 -4.9474\\
}--cycle;


\addplot[area legend, line width=1.0pt, draw=black, fill=DarkGray, forget plot]
table[row sep=crcr] {
x	y\\
9.0796 -5.0000\\
8.8754 -2.8500\\
8.4770 -0.8068\\
7.8968 1.1599\\
7.1413 3.0726\\
6.2132 4.9474\\
-9.0796 5.0000\\
-8.8754 2.8500\\
-8.4770 0.8068\\
-7.8968 -1.1599\\
-7.1413 -3.0726\\
-6.2132 -4.9474\\
}--cycle;

\end{axis}

\end{tikzpicture}

\end{minipage}\\
    (a) Tube-RMPC & (b) RMPC-ADF
    \end{tabular}
    \normalsize
	\caption{RoA of Tube-RMPC and RMPC-ADF for $t_{\mathrm{S}}=0.1$ and $N=10$ (dark gray), $N=20$ (medium gray), $N=40$ (light gray).\label{Fig:RoAComp}}
\end{figure}
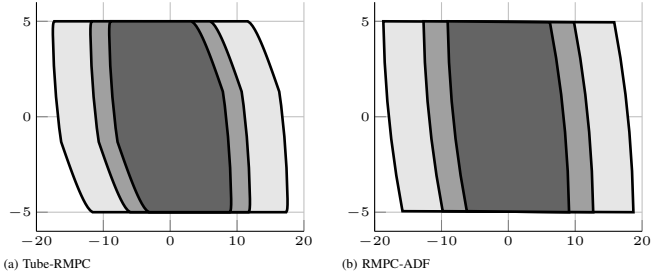

All three approaches come with the formal guarantee of input-to-state stability. The main advantage of DRMPC over Tube-RMPC and RMPC-ADF is that no computation of a (minimal or maximal) RPI set is required. Figure \ref{Fig:RoADRMPC} shows the RoA of Offline DRMPC, where the terminal set is (a) the origin and (b) the maximal PI set based on an LQR. Figure \ref{Fig:RoAComp} compares these RoAs against those of (a) Tube-RMPC and (b) RMPC-ADF.

\begin{figure}[H]
    \tiny
	\begin{tabular}{p{0.46\columnwidth}p{0.46\columnwidth}}

\begin{minipage}[t]{0.4\columnwidth}

\begin{tikzpicture}

\begin{axis}[%
width=1.400in,
height=1.200in,
at={(0.0in,0.0in)},
scale only axis,
xmin=-20,
xmax=20,
ymin=-6,
ymax=6,
label style={font=\tiny},
tick label style={font=\tiny},
axis background/.style={fill=white},
axis x line*=bottom,
axis y line*=left,
xmajorgrids,
ymajorgrids
]


\addplot[area legend, line width=1.0pt, draw=black, fill=LiteGray, forget plot]
table[row sep=crcr] {
x	y\\
13.7678 -5.0000\\
13.5635 -2.8500\\
13.1651 -0.8069\\
12.5850 1.1597\\
11.8295 3.0723\\
10.9015 4.9470\\
-13.7678 5.0000\\
-13.5635 2.8500\\
-13.1651 0.8069\\
-12.5850 -1.1597\\
-11.8295 -3.0723\\
-10.9015 -4.9470\\
}--cycle;


\addplot[area legend, line width=1.0pt, draw=black, fill=MedGray, forget plot]
table[row sep=crcr] {
x	y\\
7.5166 -5.0000\\
7.3124 -2.8500\\
6.9147 -0.8109\\
6.3366 1.1487\\
5.5849 3.0520\\
4.6623 4.9157\\
-7.5166 5.0000\\
-7.3124 2.8500\\
-6.9147 0.8109\\
-6.3366 -1.1487\\
-5.5849 -3.0520\\
-4.6623 -4.9157\\
}--cycle;


\addplot[area legend, line width=1.0pt, draw=black, fill=DarkGray, forget plot]
table[row sep=crcr] {
x	y\\
3.8783 -5.0000\\
3.6740 -2.8500\\
3.2827 -0.8434\\
2.7207 1.0617\\
1.9900 2.9117\\
1.0742 4.7617\\
0.9324 5.0000\\
-3.8783 5.0000\\
-3.6740 2.8500\\
-3.2827 0.8434\\
-2.7207 -1.0617\\
-1.9900 -2.9117\\
-1.0742 -4.7617\\
-0.9324 -5.0000\\
}--cycle;

\end{axis}

\end{tikzpicture}

\end{minipage}
    &

\begin{minipage}[t]{0.4\columnwidth}

\begin{tikzpicture}

\begin{axis}[%
width=1.400in,
height=1.200in,
at={(0.0in,0.0in)},
scale only axis,
xmin=-20,
xmax=20,
ymin=-6,
ymax=6,
label style={font=\tiny},
tick label style={font=\tiny},
axis background/.style={fill=white},
axis x line*=bottom,
axis y line*=left,
xmajorgrids,
ymajorgrids
]


\addplot[area legend, line width=1.0pt, draw=black, fill=LiteGray, forget plot]
table[row sep=crcr] {
x	y\\
18.5865 -5.0000\\
18.3823 -2.8500\\
17.9838 -0.8069\\
17.4037 1.1597\\
16.6482 3.0723\\
15.7203 4.9470\\
-18.5865 5.0000\\
-18.3823 2.8500\\
-17.9838 0.8069\\
-17.4037 -1.1597\\
-16.6482 -3.0723\\
-15.7203 -4.9470\\
}--cycle;


\addplot[area legend, line width=1.0pt, draw=black, fill=MedGray, forget plot]
table[row sep=crcr] {
x	y\\
12.3890 -5.0000\\
12.1847 -2.8500\\
11.7871 -0.8109\\
11.2090 1.1487\\
10.4572 3.0520\\
9.5347 4.9157\\
-12.3890 5.0000\\
-12.1847 2.8500\\
-11.7871 0.8109\\
-11.2090 -1.1487\\
-10.4572 -3.0520\\
-9.5347 -4.9157\\
}--cycle;


\addplot[area legend, line width=1.0pt, draw=black, fill=DarkGray, forget plot]
table[row sep=crcr] {
x	y\\
9.0518 -5.0000\\
8.8475 -2.8500\\
8.4562 -0.8434\\
7.8942 1.0617\\
7.1635 2.9117\\
6.2477 4.7617\\
6.1059 5.0000\\
-9.0518 5.0000\\
-8.8475 2.8500\\
-8.4562 0.8434\\
-7.8942 -1.0617\\
-7.1635 -2.9117\\
-6.2477 -4.7617\\
-6.1059 -5.0000\\
}--cycle;

\end{axis}

\end{tikzpicture}

\end{minipage}\\
    (a) DRMPC with terminal point & (b) DRMPC with terminal PI set
    \end{tabular}
    \normalsize
	\caption{RoA of Offline DRMPC for $t_{\mathrm{S}}=0.1$ and $N=10$ (dark gray), $N=20$ (medium gray), $N=40$ (light gray).\label{Fig:RoADRMPC}}
\end{figure}
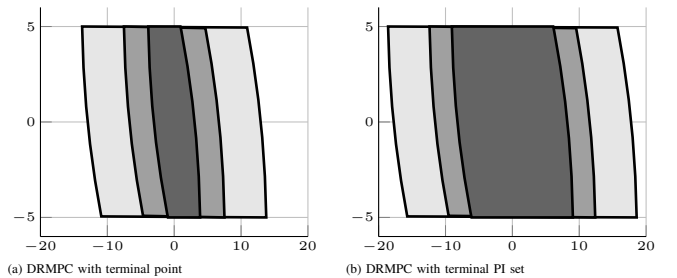

For a quantitative comparison, Tables \ref{Tab:RoAvsTubeMPC} and \ref{Tab:RoAvsADFMPC} show the excess volume of the RoA of Offline DRMPC with terminal PI set over that of Tube-RMPC and RMPC-ADF, respectively.

\renewcommand{\arraystretch}{1.4}
\begin{table}[H]
	\centering
	\begin{tabular}{|l|l||D{.}{.}{3.2}|D{.}{.}{3.2}|D{.}{.}{3.2}|D{.}{.}{3.2}|D{.}{.}{3.2}|}
		\hline
		\multicolumn{2}{|c||}{$\,$} & \multicolumn{5}{c|}{\textbf{horizon length $N$}}\\\cline{3-7}
		\multicolumn{2}{|c||}{$\,$} & \multicolumn{1}{c|}{$10$} & \multicolumn{1}{c|}{$20$} & \multicolumn{1}{c|}{$30$} & \multicolumn{1}{c|}{$40$} & \multicolumn{1}{c|}{$50$}\\\hline\hline
		\multirow{4}{0.3cm}{\rotatebox{90}{\textbf{parameter $t_{\mathrm{S}}$}}}
		 	& $0.10$ & +4.1\% & +7.6\% & +8.6\% & +8.2\% & +4.2\% \\\cline{2-7}
			& $0.15$ & +4.3\% & +5.7\% & +1.3\% & +1.3\% & +1.3\%\\\cline{2-7}
			& $0.25$ & +4.2\% & +0.6\% & +0.6\% & +0.6\% & +0.6\%\\\cline{2-7}
			& $0.40$ & +0.3\% & +0.3\% & +0.3\% & +0.3\% & +0.3\%\\\hline
	\end{tabular}
	\caption{Excess volume of RoA of Offline DRMPC with terminal PI set over Tube-RMPC\label{Tab:RoAvsTubeMPC}}
\end{table}
\renewcommand{\arraystretch}{1.0}

\renewcommand{\arraystretch}{1.4}
\begin{table}[H]
	\centering
	\begin{tabular}{|l|l||D{.}{.}{3.2}|D{.}{.}{3.2}|D{.}{.}{3.2}|D{.}{.}{3.2}|D{.}{.}{3.2}|}
		\hline
		\multicolumn{2}{|c||}{$\,$} & \multicolumn{5}{c|}{\textbf{horizon length $N$}}\\\cline{3-7}
		\multicolumn{2}{|c||}{$\,$} & \multicolumn{1}{c|}{$10$} & \multicolumn{1}{c|}{$20$} & \multicolumn{1}{c|}{$30$} & \multicolumn{1}{c|}{$40$} & \multicolumn{1}{c|}{$50$}\\\hline\hline
		\multirow{4}{0.3cm}{\rotatebox{90}{\textbf{parameter $t_{\mathrm{S}}$}}}
		 	& $0.10$ & -15.3\% & -6.1\% & -2.1\% & -0.6\% & -0.0\%\\\cline{2-7}
			& $0.15$ & -6.3\% & -1.5\% & -0.0\% & -0.0\% & -0.0\%\\\cline{2-7}
			& $0.25$ & -1.8\% & -0.0\% & -0.0\% & -0.0\% & -0.0\%\\\cline{2-7}
			& $0.40$ & -0.0\% & -0.0\% & -0.0\% & -0.0\% & -0.0\%\\\hline
	\end{tabular}
	\caption{Excess volume of RoA of Offline DRMPC with terminal PI set over ADF-MPC\label{Tab:RoAvsADFMPC}}
\end{table}
\renewcommand{\arraystretch}{1.0}

In conclusion, for the chosen example system DRMPC performs slightly better than Tube-RMPC and slightly worse than RMPC-ADF, in terms of the volume of their RoA. However, the performance of all three methods is fairly close, especially for longer prediction horizons. Thus the appropriate method may be selected based on other factors, such as implementation preferences.

\section{Conclusion}\label{Sec:Conclusion}

This paper has proposed Deadbeat Robust Model Predictive Control (DRMPC) as a new approach of Robust Model Predictive Control (RMPC) for linear systems with additive disturbances. The main idea is to completely extinguish the disturbances within a small number of time steps, called the deadbeat horizon. In two variations, the deadbeat inputs are either part of the online optimization (Online DRMPC) or pre-calculated during the design phase of the controller (Offline DRMPC). The first of more of theoretical interest, due to its computational complexity. The latter is presumably more relevant for practical applications.

Full system-theoretic guarantees have been established for both versions. The relationships of DRMPC to existing approaches, including Tube-based RMPC, RMPC with Affine Disturbance Feedback, and Minimax MPC, have been discussed and demonstrated. The performance of DRMPC turns out to be competitive with other state-of-the-art RMPC approaches. 

The main advantage of DRMPC, however, is the fact that no calculation of any Positive Invariant (PI) or Robust Positive Invariant (RPI) set is required. This implies that it is straightforward to extend the DRMPC approach to linear time-varying (LTV) and linear parameter-varying (LPV) systems. The details will be explored in future research.


\bibliographystyle{abbrv}
\bibliography{bibcontr,bibeng,bibmath}



\end{document}